\documentclass{article} 
\usepackage{iclr2021_conference,times,authblk}

\usepackage{amsmath,amsfonts,bm}

\def\eqref#1{equation~\ref{#1}}

\def\1{\bm{1}}

\DeclareMathAlphabet{\mathsfit}{\encodingdefault}{\sfdefault}{m}{sl}
\SetMathAlphabet{\mathsfit}{bold}{\encodingdefault}{\sfdefault}{bx}{n}

\newcommand{\E}{\mathbb{E}}

\newcommand{\Cov}{\mathrm{Cov}}

\usepackage{hyperref}
\usepackage{url}
\usepackage{graphicx}
\usepackage{amsmath,amssymb, array}
\usepackage[draft=true]{minted}
\usepackage{breqn}

\title{Generalized Wick Decompositions}
\author{Chris MacLeod\thanks{Author list alphabetical.}}
\author{Evgenia Nitishinskaya}
\author{Buck Shlegeris}
\affil{Redwood Research}

\newcommand{\kf}{\operatorname{K}_f}
\newcommand{\kg}{\operatorname{K}_g}
\newcommand{\wf}{\operatorname{W}_f}
\newcommand{\sbar}{[n] \setminus S}

\iclrfinalcopy 
\begin{document}

\maketitle

\section{Introduction}

The cumulant decomposition is a way of decomposing the expectation of a product of random variables (e.g. $\E[XYZ]$) into a sum of terms corresponding to partitions $\pi$ of these variables. Informally, the term corresponding to a partition measures the contribution to the expectation from irreducible interactions within each block of the partition. (That is, variables in a block of a partition interacting with each other in a way that can't be captured by a finer-grained partition.)

The Wick decomposition decomposes a product of (not necessarily random) variables into a sum of terms corresponding to subsets $S$ of the variables. Informally, each term measures the contribution from that subset of variables taking on their particular values, as compared to their expected values on a reference distribution $X_{[n]}$.

In this work, we review these two decompositions, then generalize each one to a new decomposition, where the product function is generalized to an arbitrary function $f$. 

\begin{table}[ht]
\caption{A taxonomy of interconnected decompositions to be defined further in this work.}
\label{taxonomy}
\begin{center}
\begin{tabular}{l | l}
Cumulant Decomposition & Wick Decomposition \\
$\E[\prod_{i \in [n]} X_i] = \sum_{\pi}\operatorname{K}(\pi)$ & $\prod_{i \in [n]} x_i = \sum_{S} \operatorname{W}(X_{[n]}, x_S)$ \\
\hline
Generalized Cumulant Decomposition & Generalized Wick Decomposition \\ 
$\E[f(X_1, \ldots, X_n)] = \sum_{\pi} \kf(\pi)$ & $f(x_1, \ldots, x_n) = \sum_{S} \wf(X_{[n]}, x_{[n]}, S)$ \\
\end{tabular}
\end{center}
\end{table}

\section{Expectations of Products of Random Variables}

First we'll review the expectation of a product $\E[XY]$. Rearranging the definition of covariance:
\begin{align*}
\Cov[X, Y] &= \E[(X - \E[X])(Y - \E[Y])] \\
&= \E[XY] - \E[X]\E[Y] \\
\E[XY] &= \E[X]\E[Y] + \Cov[X, Y] 
\end{align*}
Recall that if $X, Y$ are independent this implies $\Cov[X, Y] = 0$ but the converse is not true; $\Cov[X, Y] = 0$ is a weaker condition. We can think of $\Cov[X, Y] = 0$ as stating: the joint distribution doesn't have any information in it about $\E[XY]$ than isn't already in $\E[X]$ and $\E[Y]$ separately. If the covariance is positive, it means specifically that the joint distribution is denser in regions of high product than you would expect just from looking at the marginal distributions $\E[X]$ and $\E[Y]$.

Another way to put it is that the covariance is the residual after making the best possible prediction of $\E[XY]$ given only $\E[X]$ and $\E[Y]$.

This is a form of ``explanation": the reason the product has that value is because the two terms in the decomposition add to that value. For more complicated decompositions we will have many more terms, and the hope is that many terms will be zero or negligible, giving us a simple (sparse) explanation.

Next, let's generalize to the product of an arbitrary number of variables using cumulants.

\subsection{Cumulant Decomposition}

Let $\kappa(X_1, ..., X_n)$ be the $n$th cumulant; recall that $\kappa(X) = \E[X]$ and $\kappa(X, Y) = \Cov[X, Y]$. 

Let $P_n$ be the set of partitions of $n$ random variables $\{X_1, ..., X_n\}$. $P_n$ has elements $\{\pi_1, \ldots, \pi_{\lvert P_n \rvert}\}$, where each partition $\pi$ is a set of ``blocks" $\{B_1,\ldots,B_{\lvert \pi \rvert}\}$ and each block is a set of random variables. Let $\operatorname{K}(\pi) = \prod_{B\in\pi}\kappa(B)$.

For example, one member of $P_3$ is the partition $\pi = \{\{X_1\}, \{X_2,X_3\}\}$ with two blocks $B_1 = \{X_1\}$ and $B_2 = \{X_2, X_3\}$; then $\operatorname{K}(\pi) = \kappa(X_1)\kappa(X_2,X_3)$.

With this notation, we can express the expectation of products of an arbitrary number of random variables. In fact, this gives a recursive definition of cumulants in terms of lower order cumulants:
\begin{equation}\label{cumulant_def}
\operatorname \E[X_1\cdots X_n]=\sum_{\pi\in P}\operatorname{K}(\pi)=\sum_{\pi\in P}\prod_{B\in\pi}\kappa(B)
\end{equation}
For example, when $n=3$:
\begin{align*}
\E[XYZ] &= \kappa(X)\kappa(Y)\kappa(Z)\\ 
&+ \kappa(X)\kappa(Y,Z) + \kappa(Y)\kappa(X, Z) + \kappa(Z)\kappa(X, Y)\\
&+ \kappa(X, Y, Z)
\label{ntotals}
\end{align*}

\section{Generalized Cumulant Decomposition}
Next, our plan is to generalize from $\operatorname{K}(\pi)$ above to a new function $\kf(\pi)$, where $f$ is an arbitrary function of the random variables. When $f$ is the product function, $\kf$ reduces to $\operatorname{K}$. The generalized version of \eqref{cumulant_def} is:
\begin{equation}\label{gen_cumulant_def}
\operatorname \E_{(x_1,...,x_n) \sim (X_1,...,X_n)}[f(x_1,...,x_n)] := \sum_{\pi\in P}\kf(\pi)
\end{equation}
$\kf$ is defined via the base case $\kf(\{\{X\}\}) = \E_{x \sim X}[f(x)]$, plus the property that for any $B_i \in \pi$, if we let $g = \kf(\{B_i\})$ then:
\begin{equation}
\label{recur}
    \kf(\pi) := \kg(\pi \setminus B_i)
\end{equation}
\subsection{Two Variable Case}
When $f$ has two arguments, there are two partitions and thus the decomposition is:
\begin{equation}
\label{totals}
    \E_{(x,y)\sim(X,Y)}[f(x, y)] = \kf(\{\{X\}, \{Y\}\}) + \kf(\{\{X, Y\}\})
\end{equation}

Just like above, the first term represents the contribution or best estimate using the univariate distributions, and the second represents the residual or additional contribution using the joint distribution. 

Applying \eqref{recur} to $\kf(\{\{X\}, \{Y\}\})$, we let $g = \kf(\{\{X\}\})$ and then:
\begin{align*}
    \kf(\{\{X\}, \{Y\}\}) 
    &= \kg(\{\{Y\}\}) && \text{By \eqref{recur}} \\
    &= \E_{y \sim Y} [\kf(\{\{X\}\}) ] && \text{Base case on $\kg$} \\
    &= \E_{y \sim Y} [\E_{x \sim X} [f(x, y)]] && \text{Base case on $\kf$}
\end{align*}

Now we can solve \eqref{totals} for $\kf(\{\{X, Y\}\})$:
\begin{equation*}
    \kf(\{\{X, Y\}\}) = \E_{(x,y) \sim (X,Y)}[f(x, y)] - \E_{x \sim X, y \sim Y}[f(x, y)]
\end{equation*}

Just like $\Cov[X, Y] = 0$ means there's no additional information about $\E[XY]$ in the joint distribution, $\kf(\{\{X, Y\}\}) = 0$ means there's no additional information about $\E[f(X, Y)]$ in the joint distribution. Thus, we call $\kf(\{\{X, Y\}\})$ the generalized covariance.

Note that identities for the covariance need not hold for the generalized covariance. For example, $\Cov[cX,Y] = c\Cov[X,Y]$ for constant $c$, but this is due to the distributive property of the product function and won't hold for arbitrary $f$ when using generalized covariance.

\subsection{Three Arguments}

For three or more arguments, we can define $\kf$ for each partition containing more than one block by recursing on any block using \eqref{recur}. Then there's a final partition with only one block, which is defined by rearranging \eqref{gen_cumulant_def}.

We'll work the algebra here for demonstration. For brevity, we write partitions using vertical bars: $X_1|X_2,X_3$ is shorthand for $\{\{X_1\}, \{X_2,X_3\}\}$.
\begin{align}
\label{3totals}
\begin{split}
\E_{(x,y,z) \sim (X,Y,Z)}[f(x, y, z)] &= \kf(X|Y|Z) \\ 
&+ \kf(X, Y|Z) + \kf(Y, Z|X) + \kf(X, Z|Y)  \\
&+ \kf(X, Y, Z)
\end{split}
\end{align}

To evaluate $\kf(Z|X, Y)$, define $g(x, y) = \E_{z \sim Z}[f(x, y, z)]$. Then expand $\kg(X, Y)$:
\begin{align*}
\kf(X, Y|Z) &= \kg(X, Y) \\ 
&= \E_{(x,y) \sim (X,Y)}[\E_{z \sim Z}[f(x, y, z)] - \E_{x \sim X, y \sim Y}[\E_{z \sim Z}[f(x, y, z)]]
\end{align*}

Similarly:
\[ \kf(X|Y, Z) = \E_{(y,z) \sim (Y,Z)}[\E_{x \sim X}[f(x, y, z)] - \E_{y \sim Y, z \sim Z}[\E_{x \sim X}[f(x, y, z)]]  \]
\[ \kf(Y|X, Z) = \E_{(x,z) \sim (X,Z)}[\E_{y \sim Y}[f(x, y, z)] - \E_{x \sim X, z \sim Z}[\E_{y \sim Y}[f(x, y, z)]]  \]
\[\kf(X|Y|Z) = \kg(Y|Z) = \E_{x \sim X, y \sim Y, z \sim Z}[f(x, y, z)] \]

Finally, by rearranging \eqref{3totals}:
\begin{multline*}
\kf(X, Y, Z) = \\
2\E_{x \sim X, y \sim Y, z \sim Z}[f(x, y, z)] \\
- \E_{(y,z) \sim (Y,Z)}[\E_{x \sim X}[f(x, y, z)] \\
- \E_{(x,z) \sim (X,Z)}[\E_{y \sim Y}[f(x, y, z)] \\
- \E_{(x,y) \sim (X,Y)}[\E_{z \sim Z}[f(x, y, z)] \\
+ \E_{(x,y,z) \sim (X,Y,Z)}[f(x, y, z)] \\
\end{multline*}
\subsection{Matrix Form}

Each $\kf$ is a linear combination of expectations of $f$. Each expectation can be measured by sampling, and once we have all possible expectations, we can compute all possible $\kf$ terms efficiently via a change of basis matrix.

In the following we omit the arguments to $f$ and write the partition under the expectation, so $\E_{{x \sim X}, (y,z) \sim (Y,Z)}[f(x, y, z)]$ is abbreviated $\E_{X|Y,Z}[f]$. Note that each expectation uses all the variable exactly once, and the only distinction is which variables are joint and which are separate. 

\[
    \begin{bmatrix}
        \kf(X|Y) \\
        \kf(X,Y)
    \end{bmatrix} =
    \left[
        \begin{array}{cc}
            1 & 0 \\
            -1 & 1
        \end{array}
    \right]
    \begin{bmatrix}
        \E_{X|Y}[f] \\
        \E_{X,Y}[f]
    \end{bmatrix} 
\]

\[
    \begin{bmatrix}
        \kf(X|Y|Z) \\
        \kf(X|Y,Z) \\
        \kf(Y|X,Z) \\
        \kf(Z|X,Y) \\
        \kf(X|Y|Z)
    \end{bmatrix} =
    \left[
        \begin{array}{ccccc}
            1 & 0 & 0 & 0 & 0 \\
            -1 & 1 & 0 & 0 & 0 \\
            -1 & 0 & 1 & 0 & 0 \\
            -1 & 0 & 0 & 1 & 0 \\
            2 & -1 & -1 & -1 & 1
        \end{array}
    \right]
    \begin{bmatrix}
        \E_{X|Y|Z}[f] \\
        \E_{X|Y,Z}[f] \\
        \E_{Y|X,Z}[f] \\
        \E_{Z|X,Y}[f] \\
        \E_{X,Y,Z}[f]
    \end{bmatrix} 
\]

\subsection{Four or more Arguments}

Again we proceed by recursion on blocks of the partition. The only new calculation in the four argument case is $\kf(X,Y|Z,W)$, since it has two blocks of two variables each. Let $g(Z, W) = \kf(X, Y)$, then:

\begin{align*}
    \kf(X,Y|Z,W) &= \kg(Z, W) && \text{By \eqref{recur}} \\
    &= \E_{Z, W}[g(Z, W)] - \E_{Z|W}[g(Z, W)] && \text{Definition of $\kg(Z, W)$}\\
    &= \E_{Z, W}[\E_{X, Y}[f] - \E_{X|Y}[f]] - \E_{Z|W}[\E_{X, Y}[f] - [\E_{X|Y}[f]]] && \text{Expanding g} \\
    &= \E_{Z,W|X,Y}[f] - \E_{Z,W|X|Y}[f] - \E_{Z|W|X,Y}[f] + \E_{Z|W|X|Y}[f] && \text{Linearity of expectation}
\end{align*}

Again note that the order of blocks within a partition is arbitrary; we obtain the same result via $g(X, Y) = \kf(Z, W)$.

\section{Wick Product Decomposition}
Instead of decomposing the expectation of a product with cumulants, we can decompose a product of variables $x_n$ a sum of \emph{Wick terms}. Let $[n] := {\{1, \ldots, n\}}$, $X_S = (X_i)_{\{i \in S\}}$, and $x_S = (x_i)_{\{i \in S\}}$. Then each Wick term $\operatorname{W}(X_{[n]}, x_S)$ is a polynomial in $x_S$ with coefficients that are expectations of $X_{[n]}$:
\begin{equation} \label{wickdefn} \prod_{i=1}^n x_i := \sum_{S \subseteq [n]} \operatorname{W}(X_{[n]}, x_S) := \sum_{S \subseteq [n]} \E \left[\prod_{i \in [n] \setminus S} X_i\right] \varepsilon_{X_S}(x_S) \end{equation}

$\varepsilon$ is called the \emph{Wick product}; it is a polynomial in $x_S$ with coefficients that are expectations of $X_S$, defined recursively by \eqref{wickdefn} and the base case $\varepsilon_\varnothing()=1$. We can solve for $\varepsilon_{X_{[n]}}(x_{[n]})$ by subtracting off all the other terms in the sum:
\begin{equation}
\label{wickisolate}
\varepsilon_{X_{[n]}}(x_{[n]}) = \prod_{i=1}^n x_i - \sum_{S \subset [n]} \operatorname{W}(X_{[n]}, x_S) = \prod_{i=1}^n x_i - \sum_{S \subset [n]} \E \left[\prod_{i \in [n] \setminus S} X_i\right] \varepsilon_{X_S}(x_S)
\end{equation}
We often evaluate the Wick product at $(X_1,\ldots X_n)$, in which case we abbreviate the resulting random variable $\varepsilon_{X_1,\ldots,X_n}(X_1,\ldots,X_n)$ as just $\varepsilon(X_1,\ldots,X_n)$. This represents the residual after estimating the product $\Pi_{i=1}^n X_i$ using Wick products of order $<n$ and cumulants of order $<=n$.

\subsection{Basic Cases}

By \eqref{wickisolate} we have:
\begin{align*}
\varepsilon_{X_1}(x_1) = x_1 - \E[X_1]
\end{align*}
\begin{align*}
\varepsilon_{X_1, X_2}(x_1, x_2) &= x_1x_2 - \E[X_1X_2] - \E[X_2](x_1 - \E[X_1]) - \E[X_1](x_2 - \E[X_2]) \\
&= x_1x_2 - \E[X_1X_2] - \E[X_2]x_1 - \E[X_1]x_2 + 2\E[X_1]\E[X_2]
\end{align*}
\begin{align*}
\begin{split}
\varepsilon_{X_1,X_2,X_3}(x_1, x_2, x_3) &= x_1x_2x_3 \\
&- \E[X_3]x_1x_2 - \E[X_2]x_1x_3 - \E[X_1]x_2x_3  \\
&- (\E[X_2X_3] - 2\E[X_2]\E[X_3])x_1 \\
&- (\E[X_1X_3] - 2\E[X_1]\E[X_3])x_2 \\
&- (\E[X_1X_2] - 2\E[X_1]\E[X_2])x_3 \\
&- \E[X_1X_2X_3]- 6\E[X_1]\E[X_2]\E[X_3] \\
&+ 2\E[X_1]\E[X_2X_3] + 2\E[X_2]\E[X_1X_3] + 2\E[X_3]\E[X_2X_3]
\end{split}
\end{align*}

\subsection{Expectation of a Wick Product}
\label{wick_expectation_proof}
We will prove by induction that for any $n >= 1$ and any $X_{[n]}$, 
\begin{equation}\label{wick_expectation}\E[\varepsilon(X_1,\ldots,X_n)] = 0 \end{equation}

\subparagraph*{Base case:} $\E[\varepsilon(X_1)] = \E[X_1 - \E[X_1]] = 0$.
\subparagraph*{Induction step:} Assume \eqref{wick_expectation} holds for any number of variables $1 \leq n < k$. Then:
\begin{align*}
\E[\varepsilon(X_1,\ldots,X_k)] &= \E \left[ \prod_{i=1}^k X_i \right] - \sum_{S \subset [k]} \E \left[\prod_{i \in [k] \setminus S} X_i\right] \E[\varepsilon(x_S)]
\end{align*}

By the induction hypothesis, all terms in the sum vanish except the $S=\varnothing$ term: 
\begin{align*}
\E[\varepsilon(X_1,\ldots,X_k)] &= \E \left[ \prod_{i=1}^k X_i \right] - \E \left[\prod_{i \in [k] \setminus \varnothing} X_i\right] \E[\varepsilon(\varnothing)] = 0
\end{align*}
Thus, if \eqref{wick_expectation} holds for $1 \leq n <k$ then it holds for $n=k$, completing the proof.

\subsection{Derivative of a Wick Product}
We will prove by induction that for any $n \ge 1$ and any $i \in [n]$,

\begin{equation}
\label{wick_derivative}
\frac{\partial \varepsilon_{X_{[n]}}}{\partial x_i}(x_{[n]}) = \varepsilon_{X_{[n] \setminus i}}(x_{[n] \setminus i})
\end{equation}

\subparagraph*{Base case:} $\frac{\partial \varepsilon_{X_1}}{\partial x_1}(x_{1}) = \frac{\partial (x_1 - \E[X_1])}{\partial x_1} = 1 = \varepsilon_{\varnothing}$.

\subparagraph*{Induction step:} Suppose \eqref{wick_derivative} holds for any number of variables $1 \leq n < k$ and any $i \leq n$. Then:
\begin{align*}
    &{}\frac{\partial [\varepsilon(X_1,\ldots,X_k)]}{\partial x_i}(x_1, \cdots, x_k)\\
    &= \frac{\partial}{\partial x_i} \prod_{j\in [k] } x_j - \sum_{S \subset [k]} \E \left[\prod_{j \in [k] \setminus S} X_j\right]  \frac{\partial}{\partial x_i} \varepsilon_{X_S}(x_S) \\
    &= \prod_{j\in [k] \setminus i} x_j - \sum_{S \subset [k] \text{ s.t. } i \in S} \E \left[\prod_{j \in [k] \setminus S} X_j\right]  \frac{\partial}{\partial x_i} \varepsilon_{X_S}(x_S) && \frac{\partial \varepsilon}{\partial x_i} = 0 \text{ if } i \notin S \\
    &=  \prod_{j\in [k] \setminus i} x_j - \sum_{S \subset [k] \text{ s.t. } i \in S} \E \left[\prod_{j \in [k] \setminus S} X_j\right]  \varepsilon_{X_{S \setminus i}}(x_{S \setminus i}) && \text {by the I.H.} \\
    &=  \prod_{j\in [k] \setminus i} x_j - \sum_{S' \subset [k] \setminus i} \E \left[\prod_{j \in ([k] \setminus i) \setminus S'} X_j\right]  \varepsilon_{X_{S'}}(x_{S'}) && \text {Let } S'=S \setminus i \\
    &= \varepsilon_{X_{[k] \setminus i}}(x_{[k] \setminus i})
\end{align*}
\section{Generalized Wick Decomposition}

Generalized Wick products are the bottom right corner in our taxonomy (Figure \ref{taxonomy}). As such they can be viewed in two ways:
\begin{itemize}
\item Generalized Wick products are to Wick products as generalized cumulants are to cumulants: both move from explaining products to explaining arbitrary functions.
\item Generalized Wick products are to generalized cumulants as Wick products are to cumulants: both move from explaining expectations of random variables $X_n$ to explaining individual values $x_n$.
\end{itemize}

We can decompose $f$ into a sum of \emph{generalized Wick terms} $\wf(X_{[n]}, x_{[n]}, S)$. A regular Wick term was able to depend only on $x_S$, but for an arbitrary function $f$ we need a complete set of arguments $x_{[n]}$ in order to evaluate $f$. Now the purpose of $S$ is to determine which arguments of $f$ use the individual values $x_i$ and which arguments use draws from the expectation $d_i \sim X_i$. The decomposition is: 
\begin{equation}
\label{gwickdefn}
f(x_1,\ldots,x_n) := \sum_{S \subseteq [n]} \wf(X_{[n]}, x_n, S) := \sum_{S \subseteq [n]} \E_{d_{\sbar} \sim X_{\sbar}} \left[ \omega_{f, X_S} \right] 
\end{equation}
We call $\omega_{f}$ the \emph{generalized Wick product}, defined recursively by \eqref{gwickdefn} and the base case $\omega_{f, \varnothing} = f$. While a Wick product contains products of $x_s$, the generalized version is not a literal product. For a fixed $X_S$, the expression $\omega_{f, X_S}$ is a function with the same type as $f$. As before, we can isolate $\omega_f$ by subtracting off all other terms:
\begin{equation}
\label{genwickisolate}
\omega_{f, X_{[n]}}(x_1,\ldots,x_n) = f(x_1,\ldots,x_n) - \sum_{S \subset [n]} \E_{d_{\sbar} \sim X_{\sbar}} \left[ \omega_{f, X_S} \right]
\end{equation}
\subsection{Basic Cases}
We'll write $\E_{\pi}$ for an expectation over $\pi$ as before, but now $\pi$ isn't necessarily a complete partition, in which case $\E_{\pi}$ is a function with free variables to be inferred from context. For example, the expression $\E_{d_1 \sim X_1, (d_2, d_3) \sim (X_2, X_3)} f(d_1, d_2, d_3, x_4, x_5)$ is abbreviated $\E_{1|23}$. 

Note that the equations are identical to the basic Wick product, except that an expression in the basic product like $\E_{1|23}[x_4x_5]$ translates to $\E_{1|23}$ (where $x_4$ and $x_5$ are free variables). For example, for $f(x_1)$:
\begin{align*}
\omega_{f, X_1}(x_1) &= f(x_1) - \E_{d_1 \sim X_1} [f(d_1)] = f(x_1) - \E_1\\
\end{align*}
For $f(x_1, x_2)$:
\begin{align*}
\omega_{f, X_1}(x_1, x_2) &= f(x_1, x_2) - \E_{d_1 \sim X_1} [f(d_1, x_2)] = f(x_1, x_2) - \E_1\\
\end{align*}
\begin{dmath*}
{} \omega_{f, X_1, X_2}(x_1, x_2)
= f(x_1, x_2) - \E_{d_2 \sim X_2}[\omega_{f, X_1}] - \E_{d_1 \sim X_1}[\omega_{f, X_2}] - \E_{(d_1, d_2) \sim (X_1, X_2)}[f(d_1, d_2)]
= f(x_1, x_2) -\E_{d_2 \sim X_2}[f(x_1, d_2)] - \E_{d_1 \sim X_1}[f(d_1, x_2)] + 2\E_{d_1 \sim X_1, d_2 \sim X_2}[f(d_1, d_2)] - \E_{(d_1, d_2) \sim (X_1, X_2)}[f(d_1, d_2)]
= f(x_1, x_2) - \E_2 - \E_1 + 2\E_{1|2} - \E_{12}
\end{dmath*}
\subsection{Properties}
For any $f$ and any $n >=1$,
\begin{equation}
\E_{(d_1,\ldots,d_n) \sim (X_1, \ldots, X_n)}[\omega_{f, X_{[n]}}(d_1,\ldots,d_n)] = 0
\end{equation}
The proof by induction is identical to the Wick product case in Section \ref{wick_expectation_proof}.

Note that we could also define
\begin{equation}
\wf(X_{[n]}, x_{[n]}, S, \pi_{\sbar}) := K_{\omega_{f, X_S}} (\pi_{\sbar})
\end{equation}
so the usual $\wf$ is a shorthand for the trivial partition case $\wf(X_{[n]}, x_{[n]}, S, \{\sbar\})$.

\subsubsection{Commutativity of expectation and generalized Wick product}
In \eqref{gwickdefn} we defined $\wf$ as an expectation over generalized Wick products. In fact, this is equal to a generalized Wick product of a marginalized $f$:
\begin{equation}
\E_{d_T \sim X_T} \left[ \omega_{f, X_S} \right] = \omega_{g, X_S}
\end{equation}
where $g = \E_{d_T \sim X_T} f$ is a function from $X_{[n] \setminus T}$ and $S \cap T = \varnothing$

\subparagraph*{Base case:}
\begin{align*}
\E_{d_T \sim X_T} \left[ \omega_{f, X_\varnothing} \right] &= \E_{d_T \sim X_T} f = g \\
&= \omega_{g, \varnothing}
\end{align*}

\subparagraph*{Induction step:}
By induction, the claim holds for all $U \subset S$.
\begin{align*}
\E_{d_T \sim X_T} \left[ \omega_{f, X_S} \right] &= \E_{d_T \sim X_T} \left[ f(x_1,\ldots,x_n) - \sum_{U \subset S} \E_{d_{S \setminus U} \sim X_{S \setminus U}} \left[ \omega_{f, X_U} \right] \right]\\
&= g - \sum_{U \subset S} \E_{d_{S \setminus U} \sim X_{S \setminus U}} \E_{d_T \sim X_T} \omega_{f, X_U} \\
&= g - \sum_{U \subset S} \E_{d_{S \setminus U} \sim X_{S \setminus U}} \omega_{g, X_U}\\
&= \omega_{g, X_S}
\end{align*}

\subsubsection{Derivative of a Generalized Wick Product}

We will prove by induction that if $f$ is differentiable, for any $n >= 1$ and any $X_{[n]}$ satisfying suitable regularity conditions,

\begin{equation}
\frac{\partial}{\partial x_i} \omega_{f, X_{[n]}}(x_{[n]}) = \omega_{\frac{\partial f}{\partial x_i}, X_{[n]}}(x_{[n]})
\end{equation}

\subparagraph*{Base case:}
\begin{equation*}
\frac{\partial}{\partial x_1} \omega_{f, X_1}(x_1) = \frac{\partial f}{\partial x_1}(x_1) - \frac{\partial}{\partial x_1} \E_{d_1 \sim X_1} [f(d_1)]]
\end{equation*}
\begin{equation*}
\omega_{\frac{\partial f}{\partial x_1}, X_1}(x_1) = \frac{\partial f}{\partial x_i}(x_1) - \E_{d_1 \sim X_1} [\frac{\partial f}{\partial x_1}(d_1)]
\end{equation*}
\subparagraph*{Induction step:}
\begin{align*}
&{} \frac{\partial}{\partial x_i} \omega_{f, X_{[k]}}(x_{[k]}) \\
&= \frac{\partial f}{\partial x_i} - \sum_{S \subset [k]} \frac{\partial}{\partial x_i} \E_{x_{\sbar} \sim X_{\sbar}} \left[ \omega_{f, X_S} \right] \\
&= \frac{\partial f}{\partial x_i} - \sum_{S \subset [k]} \E_{x_{\sbar} \sim X_{\sbar}} \left[ \frac{\partial}{\partial x_i}\omega_{f, X_S} \right] && \text{Interchange derivative and expectation} \\
&= \frac{\partial f}{\partial x_i} - \sum_{S \subset [k]} \E_{x_{\sbar} \sim X_{\sbar}} \left[ \omega_{\frac{\partial f}{\partial x_i}, X_S}(x_S) \right] && \text{I.H.} \\
&= \omega_{\frac{\partial f}{\partial x_i}, X_{[k]}}(x_{[k]})
\end{align*}


\section{Conclusion}

We've defined a function $\kf$ which acts like a generalized ``product of cumulants" function. Given a partition with some number of blocks, it tells you the contribution to $\E[f]$ when members of a block are drawn jointly, over and above the contributions of lower-order expectations.

You can compute all the $\kf$ via a linear transformation of all the different ways to sample expectations of f. 

An extension of this can be applied to resampling on a tree-ified graph. Say that X, Y, and Z are input nodes in the tree-ified graph; then we can define some intuitive measure of how ``important" it was that X and Y (corresponding to two particular paths) were sampled together, over and above their individual importance to the graph's output.

$\wf$ does a similar thing, but accounts for not just the distributions $X_i$ but particular subsets of values $x_S$. This allows us to do ``path patching attribution" \cite{pathpatching} where we measure how important that specific combination of inputs was to producing the output.

\subsubsection*{Acknowledgements}

Thanks to Joe Benton and Ryan Greenblatt for previous work on Wick products (unpublished) and some of the properties they proved, which were helpful for reference. Thanks to Nicholas Goldowsky-Dill, Jacob Hilton, and Fabien Roger for feedback.

\bibliography{References}
\bibliographystyle{plainnat}

\section{Appendix - Python implementation}

Full listing including tests at: \url{https://github.com/redwoodresearch/cumulant_decomposition}

\begin{minted}[
    fontsize=\footnotesize,
    frame=lines
]{python}

from pprint import pprint
from typing import cast, TypeVar, Sequence
from collections import defaultdict
from functools import lru_cache

import numpy as np
from sympy.utilities.iterables import multiset_partitions

Block = frozenset[int]
Partition = frozenset[Block]
T = TypeVar("T")


def part(p: list[list[int]]) -> Partition:
    """Convert p to a hashable representation of a partition."""
    return Partition([Block(b) for b in p])


@lru_cache(maxsize=None)
def partitions(items: frozenset[int]) -> list[Partition]:
    """Return all partitions of items."""
    parts = cast(list[list[list[int]]], multiset_partitions(list(items)))
    return [Partition([Block(b) for b in part]) for part in parts]


@lru_cache(maxsize=None)
def factorial(n: int) -> int:
    if n == 0 or n == 1:
        return 1
    return n * factorial(n - 1)


def remove_zeros(d: dict[T, int]) -> dict[T, int]:
    return {k: v for k, v in d.items() if v != 0}


@lru_cache(maxsize=None)
def kf(pi: Partition) -> dict[Partition, int]:
    first, *rest = list(pi)
    if rest:
        out: defaultdict[Partition, int] = defaultdict(int)
        # note: could recurse on any block
        for b0, coef0 in kf(frozenset([first])).items():  
            for b1, coef1 in kf(frozenset(rest)).items():
                # E_b0|b1 = coef0 * E_b0 [ coef1 * E_b1[f]]
                out[frozenset([*b0, *b1])] += coef0 * coef1  
        return remove_zeros(out)

    if len(first) == 1:
        # Trivial case of K_f({{X}}) = E_X
        return {frozenset([first]): 1}  

    out: defaultdict[Partition, int] = defaultdict(int)
    for other in partitions(first):
        if pi == other:
            out[pi] = 1
        else:
            for ob, on in kf(other).items():
                out[ob] -= on  # K_f[all] = E_all - (other K_fs)
    return remove_zeros(out)

Expectation = frozenset[frozenset[int]]
Xs = frozenset[int]
Term = tuple[Expectation, Xs]

def mul_expectation(e1: Expectation, e2: Expectation) -> Expectation:
    return frozenset(term for term in chain(e1, e2) if term)

def wick(x_s: Block) -> dict[Term, int]:
    """Return a Wick decomposition such that the terms sum to the product of x_s."""
    if not x_s:
        return {(frozenset([frozenset()]), frozenset()): 1}

    terms: dict[Term, int] = defaultdict(int, {(frozenset([frozenset()]), x_s): 1})
    for subset in powerset(x_s, min=0, max=len(x_s) - 1):  # strict subsets
        for (sub_ex, sub_xs), sub_coef in wick(subset).items():
            term = mul_expectation(frozenset([x_s - subset]), sub_ex), sub_xs
            terms[term] -= sub_coef
    return terms

\end{minted}

\section{Appendix - Connection between generalized cumulants and path patching}

\subsection{Notational definitions}
\newcommand{\refine}{\mathcal{R}}

\begin{itemize}
    \item Let $P(S)$ represent the partitions of the set $S$.
    \item If $B = \{x_1, x_2, \cdots, x_n\}$ then let $\E_B f$ mean $\E_{(x_1, \cdots, x_n) \sim (X_1, \cdots, X_n)}f$. That is, the expectation over all variables in $B$ sampled together.
    \item For a set of blocks $S=\{B_1, B_2, \cdots, B_n \}$ we will write $\E_S f$ as shorthand for $\E_{B_1}\E_{B_2}\cdots \E_{B_n} f $. Intuitively, we are sampling variables together within each block. If not all variables are present, this represents a partial function.
    \item Let $\refine(\pi)$ be the set of \href{https://en.wikipedia.org/wiki/Partition_of_a_set#Refinement_of_partitions}{refinements} of a partition $\pi$. (A refinement partition may sub-partition every block of $\pi$).
\end{itemize}

\subsection{Lemma}
Let $\pi$ be a partition of the arguments of $f$, divided into non-overlapping sets of blocks $S_1$, $S_2$. We claim:
\[K_{\E_{S_1}f}(S_2) = \E_{S_1}K_{f}(S_2) \]

To see this, let $g = \E_{S_1}f$. Then
\begin{align*}
    K_{\E_{S_1}f}(S_2) &= \kg(S_2) \\
\intertext{There exists some closed form expression for this, of the form:}
   &= \sum_i c_i \E_{\alpha_i} g \\
\intertext{ where each $c_i$ is a constant, and each $\alpha_i$ partitions the variables of $S_2$. Then, substituting $g$ we get:}
   &= \sum_i c_i \E_{\alpha_i} \E_{S_1}f \\
\intertext{We can distribute out this inner expectation, as $\alpha_i$ and $S_i$ contain disjoint arguments of $f$:}
   &= \E_{S_1} \sum_i c_i \E_{\alpha_i} f \\
   &= \E_{S_1} \kf(S_2)
\end{align*}
finishing our proof.

\subsection{Theorem}
Let $\pi=\{B_1, B_2, \cdots B_n \}$ be a partition of the arguments of $f$. We claim
\[
    \E_\pi f = \sum_{\alpha \in \refine(\pi)} \kf(\alpha)
\]
 We prove this by induction.

The \textbf{base case} is immediate from the $\kf$ decomposition of the expectation. Our partition $\pi$ is then a single block $B_1 = \{x_1 \cdots x_n \}$. Note $\refine(\pi) = P(B_1)$, and so
\[
    \E_{B_1} f = \E_{x_1 \cdots x_n \sim (X_1 \cdots X_n)} f = \sum_{\alpha \in P(B_1)} \kf(\alpha) = \sum_{\alpha \in \refine(\pi)} \kf(\alpha)
\]

For the \textbf{inductive case}, let us define $S_2 = \{B_2, \cdots, B_n\}$. Then:
\begin{align*}
    \E_{\pi} f &= \E_{B_1}\E_{S_2} f \\
    &= \sum_{\alpha_1 \in P(B_1)} K_{\E_{S_2}f}(\alpha_1) \\
    &= \sum_{\alpha_1 \in P(B_1)} \E_{S_2} K_{f}(\alpha_1) && \text{by the lemma} \\
    &= \sum_{\alpha_1 \in P(B_1)} \sum_{\alpha_2 \in \refine(S_2)} K_{\kf(\alpha_1)}(\alpha_2) && \text{by inductive case} \\
    &= \sum_{\alpha_1 \in P(B_1)} \sum_{\alpha_2 \in \refine(S_2)} \kf(\alpha_1 | \alpha_2) && \text{by (\ref{recur})} \\
\intertext{However, the set of $\alpha_1 | \alpha_2 $ for all $ (\alpha_1, \alpha_2) \in P(B_1) \times \refine(S_2)$ is exactly $\refine(\pi)$ and so:}
    &= \sum_{\alpha \in \refine(\pi)} \kf(\alpha) \\
\end{align*}

\subsection{Application to Path Patching}

Let $f(x_1, \cdots, x_n, y)$ be a treeified network. Let $S \subset \{x_i\}$ be the paths we hypothesize are important, with $y \in S$ where $y$ is the label that determines the loss.

In path patching, we test to see if:
\[ \E f \stackrel{?}{\approx} \E_S \E_{\sbar} f = \E_{\{S, \sbar \}} \]

By the above theorem, this is equal to:
\[   \sum_{\alpha \in \refine(\{S, \sbar \})} \kf(\alpha)   \]

Intuitively, this corresponds to the claim that you can estimate $\E f$ well by only considering interactions between the variables in $S$ and between the variables in $\sbar$. Perhaps most importantly, this does not allow interactions between the inputs in $S$ and $y$.

This could be extended to express variations on path patching, where all unimportant inputs are not sampled coherently, by subdividing $\sbar$ into further blocks.

\end{document}